\documentclass[11pt]{article}
\usepackage{amssymb}
\usepackage{amsfonts}
\usepackage{amsmath}
\usepackage{latexsym}
\usepackage{epsfig}
\newtheorem {theorem}{Theorem}[section]
\newtheorem{open}{Open Problem}
\newtheorem{lemma}[theorem]{Lemma}
\newtheorem{proposition}[theorem]{Proposition}

\newtheorem{fact}[theorem]{Fact}
\newtheorem{definition}[theorem]{Definition}
\newtheorem{remark}[theorem]{Remark}

\newenvironment{proof}{\noindent \textbf{Proof:}}{$\Box$}
\newcommand{\bits}{\{0,1\}}
\newcommand{\fits}{\{+1,-1\}}
\newcommand{\maj}[1]{\mathsf{MAJ}_{#1}}
\renewcommand{\Pr}{{\bf P}}
\newcommand{\E}{{\bf E}}
\renewcommand{\P}{{\bf P}}

\newcommand{\Ex}{\mathop{\bf E}}
\newcommand{\p}{\mathcal{P}}
\newcommand{\A}{\mathcal{A}}
\newcommand{\B}{\mathcal{B}}
\newcommand{\C}{\mathcal{C}}

\def\R{{\Bbb R}}

\def\Z{{\Bbb Z}}
\newcommand{\eps}{\epsilon}

\begin{document}
\title{Coin flipping from a cosmic source: On error correction of truly random bits}
\author{Elchanan Mossel\thanks{The work reported here was conducted 
while the author was a PostDoc at the Theory Group, Microsoft Research} \\ Statistics, U.C. Berkeley \\ mossel@stat.berkeley.edu \and Ryan
O'Donnell\thanks{Supported by NSF grant CCR-99-12342.}\\ MIT Department of Mathematics\\
odonnell@theory.lcs.mit.edu}
\date{\today}
\maketitle
\begin{abstract}
We study a problem related to coin flipping, coding theory,
and noise sensitivity.
Consider a source of truly random bits $x \in \bits^n$, and $k$ parties,
who have noisy versions of the source bits $y^i \in \bits^n$, where for all $i$ and $j$,
it holds that $\Pr[y^i_j = x_j] = 1 - \eps$, independently for all $i$ and $j$.  That is, each party sees each bit correctly with probability $1-\epsilon$, and incorrectly (flipped) with probability
$\epsilon$, independently for all bits and all parties.
The parties, who cannot communicate, wish to agree beforehand on {\em balanced} functions $f_i : \bits^n \to \bits$ such that
$\Pr[f_1(y^1) = \ldots =  f_k(y^k)]$ is maximized. In other words, each party wants to toss a fair coin so that the
probability that all parties have the same coin is maximized. The functions $f_i$ may be thought of as an error correcting
procedure for the source $x$.

When $k=2,3$ no error correction is possible, as the optimal protocol is given by $f_i(x^i) = y^i_1$.
On the other hand, for large values of $k$, better protocols exist.
We study general properties of the optimal protocols and the asymptotic behavior of the problem with
respect to $k$, $n$ and $\eps$.
Our analysis uses tools from probability, discrete Fourier analysis, convexity and discrete symmetrization.
\end{abstract}

\section{Introduction}

Consider a source of truly random bits $x \in \bits^n$, which is accessible to $k$ parties.
If the $k$ parties want to use the source in order to obtain a common single random bit, they can easily
do so by deciding beforehand to let the common bit be $x_1$. More generally, they can decide beforehand
on any balanced function $f : \bits^n \to \bits$, and let the common bit be $f(x)$.
We call a function $f$ {\em balanced} if $\Pr_x [f(x) = 0] = \Pr_x[f(x) = 1] = 1/2$.

In this setting, there is no real advantage in taking the function $f$ to be anything other than $f(x) = x_1$.
The problem becomes more interesting when the parties receive {\em noisy} versions of the random bits.
That is, party $i$ receives $y^i$, where the bits of $y^i$ satisfy $\Pr[y^i_j = x_j] = 1 - \eps$,
independently for all $i$ and $j$. We also assume that the parties cannot communicate.
Yet, the parties want to toss the same fair coin given their noisy versions of the source.
We will now allow each party $i$ to use a different {\em balanced} function $f_i : \bits^n \to \bits$
as a coin-tossing procedure.  We want to maximize $\P[f_1(y^1) = \ldots = f_k(y^k)]$.

This problem is motivated naturally by several models in cryptography.
Think of a long one-time pad which is distributed to parties with a small probability of error.
The parties still want to use this one-time pad as their key to an encryption algorithm,
by dividing the one-time pad into blocks of length $n$,
and applying some function on each block to obtain a shorter one-time pad which has high probability
of being the same for all parties.
One setting in which such a procedure should be useful is Ding and Rabin's ``everlasting security'' \cite{DR01},
a strong encryption algorithm in the bounded storage model.
This model presupposes the existence of a satellite broadcasting a continuous stream of a huge number of random bits.
It is natural to expect some error in any reception of this stream.
A somewhat related cryptographic problem was studied by Maurer \cite{M97};

Our problem is also of interest as a noncryptographic collective
coin flipping problem.  One example of such a problem is the full
information model, introduced by Ben-Or and Linial \cite{BL90} and
studied extensively (see, e.g., the survey \cite{D00}).  In this
problem, many parties try to agree on a single random bit; each
generates a random coin toss, and there is a single protocol
(function) taking all the coin tosses and producing a bit.  The
difficulty arises from the assumption that some parties are corrupt
and can choose their coins adversarially. In our problem, the
major difference is that the parties do not communicate any random
bits, so they each must apply a protocol to a shared random
string.  And, instead of arbitrary corruptions, we assume random
ones.

The question presented in this paper is also a natural
question regarding error correcting for the broadcast channel (see e.g. \cite{CT91}) with a truly random source.
Naturally, when the source is truly random, error correction is impossible.
However, here instead of requiring that all parties receive the information transmitted to them with high probability,
we require that all parties attain the same information with high probability, and that this mutual information has high entropy.

Finally, a basic motivation comes from the study of noise-sensitivity, see \cite{KKL88} and \cite{BKS99}.
The functions $f_i$ that maximize the probability $\P[f_1(y^1) = \ldots = f_k(y^k)]$ are
in an intuitive sense stable to noise, and it turns out that when the number of parties $k$ is 2 or 3,
this intuition can be used in order to prove that the optimal
functions are just the first-bit function. Our results
for larger values of $k$ are part of an initial attempt for bounding
high norms of the Bonami-Beckner operators in the range where the
results of \cite{Bo70} and \cite{Be75} do not apply.

\subsection{Definitions and notation}
\begin{definition}\
\begin{itemize}
\item
{\bf The model}
Let $k \geq 1$ be the number of {\em parties},
and $n \geq 1$ be the {\em block length}.
Let $\eps \in (0,1/2)$ be the {\em corruption probability}.
Our space is the space of all sequences $(x,y^1,\ldots,y^k) \in \bits^{n \times (k+1)}$,
where $x$ represents the {\em source} and is chosen uniformly at random from $\bits^n$.
For each $i$, $y^i$ represents the bits that party $i$ holds and it is assumed that for all $1 \leq i \leq k$ and
$1 \leq j \leq n$, it holds that $\P[y^i_j = x_j] = 1 - \eps$, independently for all $i$ and $j$. This is our
probability space, and when we write $\P$ (or $\E$) we mean the probability (expected value) in this space.
\item
{\bf Balanced and antisymmetric functions}
Let $\B_n$ denote the set of {\em balanced }functions $f : \bits^n \to \bits$; i.e., those with $|f^{-1}(0)| = |f^{-1}(1)|$.
Let $\A_n$ denote the set of {\em antisymmetric} $n$-bit
boolean functions; i.e., those satisfying $f(\bar{x}_1
\bar{x}_2 \cdots \bar{x}_n) = \overline{f(x)}$, where the bar
denotes flipping 0's and 1's, so $\bar{x} = 1 - x$.
Note that $\A_n \subset \B_n$.
\item
{\bf Protocols}
A {\em protocol} consists of $k$ functions $f_i \in \B_n$.
An {\em antisymmetric protocol} consists of $k$ functions $f_i \in \A_n$.
For a protocol $(f_1,\ldots,f_k)$, we write $\p(f_1, \ldots, f_k ; \eps)$ for the
probability that all functions agree, so
\[
\p(f_1, \ldots, f_k ; \eps) = \P[f_1(y^1) = \ldots = f_k(y^k)].
\]
We write $\p_k(f ; \eps)$, in place of $\p(f_1, \ldots, f_k ; \eps)$, if $f_1 = \ldots = f_k = f$.
\end{itemize}
\end{definition}

It turns out that restricting all $f_i$ to be balanced is neither
necessary nor sufficient for ensuring that the output bit, when
agreed upon, is uniformly random --- see
Proposition~\ref{prop:unbalanced} for a counterexample to sufficiency.  A
sufficient condition is that every function be antisymmetric, since if all the
function are antisymmetric, then
\[
\P[f_1(y^1) = \ldots = f_k(y^k) = 1] =
\P[f_1(\overline{y^1}) = f_k(\overline{y^k}) = 0] = \P[f(y^1) = \ldots = f(y^k) = 0],
\]
where the first equality follows from the fact that $f_i$ are antisymmetric and the second
since $\P$ assigns the same probability to $(x,y^1,\ldots,y^k)$
as it does to $(\bar{x},\overline{y^1},\ldots,\overline{y^k})$. We are not aware of a weaker condition than antisymmetry that ensures that the output bit when agreed upon is uniformly random.\\

We end this section with a few more definitions.  For $S \subseteq [n]$ and $\pi$ a permutation of $[n]$, let $\pi_S : \bits^n \to \bits^n$ be defined by $\pi_S(x)_i = x_{\pi(i)}$ if $i \in S$, and $\pi_S(x)_i = \overline{x_{\pi(i)}}$ if $i \notin S$.  Any $\pi_S$ merely permutes coordinates, and flips the roles and $0$ and $1$ on some coordinates.  It's therefore easy to see that $\p(f_1 \circ \pi_S, \ldots, f_k \circ \pi_S ; \eps) = \p(f_1, \ldots, f_k ; \eps)$ for any $\pi_S$.

In order to express uniqueness results cleanly, we abuse language in the following way:  For particular $k$, $n$, and $\eps$, we say that $(f_1, \ldots, f_k)$ is the unique best protocol ``up to $\pi_S$'' if the set of best protocols is exactly $\{(f_1 \circ \pi_S, \ldots, f_k \circ \pi_S) : S \subseteq [n], \pi \in S_n\}$.

\begin{remark}
Those familiar with the Bonami-Beckner operator
(see Definition \ref{def:beckner}) will note that
$\p_k(f;\eps) = \| T_{\eps} f \|_k^k + \| T_{\eps} (1 - f)\|_k^k$.
Therefore we are looking for balanced boolean functions having large $k$'th
norm for various values of $k$.
\end{remark}
\subsection{Main results}

Methods of discrete Fourier analysis (see \cite{KKL88,BKS99,MO02} for background) give an exact solution to our problem in the cases $k=2,3$, and the best protocol, up to $\pi_S$,  is for all parties to use the function $f(x) = x_1$.  We attribute the case $k=2$ in the following theorem to folklore.

\begin{theorem} \label{thm:fourier}
For all $k$, $n$, $\eps$, if we wish to maximize the expression
\begin{equation} \label{eq:nagreements}
\E[\#(i,j) : f_i(y^i) = f_j(y^j)],
\end{equation}
the unique best protocol up to $\pi_S$ is given by $f_1 = \ldots = f_k = f$, where $f(x) = x_1$.
In particular, if $k=2$ or $k=3$, then for all $n$ and $\eps$, the unique best protocol up to $\pi_S$ for maximizing $\p(f_1,\ldots,f_k;\eps)$ is given by $f_1 = \ldots = f_k = f$, where $f(x) = x_1$.
\end{theorem}
In general we do not know how to find the optimal protocol. However, we can prove some general properties of the protocols which maximize $\p(f_1,\ldots,f_k;\eps)$. Recall that a function $f$ is {\em monotone} if for all $x, y \in \bits^n$, we have $f(x) \leq f(y)$ whenever $x \preceq y$ (in the sense $x_i \leq y_i$ for all $i$). For $x,y \in \bits^n$, we write $x \preceq_L y$ if $\sum_{i=1}^m x_i \leq \sum_{i=1}^m y_i$ for every $m=1\dots n$.
We call a function $f$ {\em left-monotone} if $f(x) \leq f(y)$ whenever $x \preceq_L y$. Note that the partial order induced by $\preceq_L$ is a refinement of the partial order induced by $\preceq$; in particular, every
left-monotone function is monotone.

The following theorem is based on convexity and on the Steiner symmetrization principle (see e.g. \cite{T93} for background).

\begin{theorem} \label{thm:general}
For all $k$, $n$, and $\eps$, any protocol which maximizes $\p(f_1,\ldots,f_k;\eps)$ among all protocols satisfies
$f_1 = \ldots = f_k = f$, where $f$ is left-monotone (up to $\pi_S$).  This theorem remains true if the phrase ``protocol'' is everywhere replaced by ``antisymmetric protocol''.
\end{theorem}

So far we haven't ruled out the possibility that the optimal protocol always consists of taking just one bit.  For $r$ an odd number, let $\maj{r}$ denote the majority function on the first $r$ bits; i.e., $\maj{r}(x)$ is 1 if $\sum_{i=1}^r x_i > r/2$ and $\maj{r}(x) = 0$ if $\sum_{i=1}^r x_i < r/2$.  Using a coupling argument, we prove
the following result:
\begin{theorem} \label{thm:majn}
For all $n$ odd, and all $\eps$, there exists a $K = K(n,\eps)$ such that for $k \geq K$,
the unique best protocol up to $\pi_S$ is given by
$f_1 = \ldots = f_k = \maj{n}$.  Moreover, as $k \to \infty$,
\begin{equation} \label{eq:as_maj}
\p_k(\maj{n};\eps) = \Theta \left( \left( 1 - \P[Bin(n,\eps) > n/2] \right)^k \right),
\end{equation}
where $Bin(n,\eps)$ is a binomial variable with parameters $n$ and $\eps$.
(This should be compared to $\Theta((1 - \eps)^k)$ for the function $f(x) = x_1$.)  When $n$ is even, a similar result is true; in place of $\maj{n}$, one should take any balanced function $f$ which has $f(x) = 1$ whenever $|\{ i : x_i = 1\}| > |\{i : x_i = 0\}|$, and $f(x) = 0$ whenever $|\{ i : x_i = 0\}| > |\{i : x_i = 1\}|$.
\end{theorem}

A dual result is obtained by fixing $n$ and $k$, and letting $\eps$ be either close to $0$ or close to $1/2$.
\begin{theorem} \label{thm:maj1}
For all $k$ and $n$, there exist $0 < \eps' = \eps'(n,k) < \eps'' = \eps''(n,k) < 1/2$, such that for all $0 < \eps < \eps'$, or $\eps'' < \eps < 1/2$, the unique best protocol up to $\pi_S$ is given by $f_1 = \ldots = f_k = f$, where $f(x) = x_1$; i.e., $f = \maj{1}$.
\end{theorem}

It may now seem like the optimal protocol consists of either taking all functions to be $\maj{n}$ or all functions to be $\maj{1}$. This is not the case however, as a computer-assisted proof shows that sometimes $\maj{r}$ is better than $\maj{1}$ and $\maj{n}$ for $1 < r < n$. See Proposition~\ref{prop:majr}.

Despite Theorem~\ref{thm:majn}, it is not true that as $k \to
\infty$, the success probability of the best protocol goes to 0
exponentially fast in $k$ (treating $\epsilon$ as fixed).  In
fact, if we allow $n$ to be an unbounded function of $k$, then the best protocol's success probability is at least inverse-polynomially large in $k$.

\begin{theorem} \label{thm:inverse-poly}
Fix $\eps$.  Then there exists a sequence $(n_k)$ such that
\[
\p_k(\maj{n_k} ; \eps) \geq \Omega\left(\frac{1}{k^{2.01/(1-2\eps)^2}}\right).
\]
It suffices to have $n_k = O(k^{4.01/(1-2\eps)^2})$.
\end{theorem}

Finally, it is natural to ask if the optimal function is always $\maj{r}$ for some $1 \leq r \leq n$ (assuming, say, $n$ is odd).\\

\noindent {\bf Conjecture M:} \emph{For a particular $k$, $\eps$, and odd $n$, there is a function $f \in \A_n$ which is
not a majority function such that $\p_k(f; \eps) > \p_k(\maj{m} ; \eps)$ for all odd $m \leq n$.}

\noindent {\bf Conjecture O:} \emph{For any given $k$, $\eps$, and odd $n$, there is an odd $m \leq n$
such that the best antisymmetric function for the parties is $\maj{m}$.}\\

In fact, we know of no counterexample to Conjecture O even if we allow the parties to use any balanced function (which could allow for a biased output).
Some evidence that resolving this conjecture could possibly be hard:
One, it is not true that for any non-majority function $f$,
and any fixed $k$, there is a majority function which dominates $f$ over all $\eps$ ---
we have a computer-verified counterexample.
Two, for certain $k, \eps$, $\p_k(\maj{n}, \eps)$ is not even unimodal as a function of $n$.
E.g., for $k = 12$, $\eps = 0.1$, the success probability decreases between $\maj{1}$ and $\maj{3}$,
increases up to $\maj{11}$, and then decreases again out to $\maj{17}$
(and appears to continue decreasing from this point on).

We conclude the introduction with a road map to the following sections.
In Section \ref{sec:fourier} we prove Theorem \ref{thm:fourier} using Fourier analysis.
In Section \ref{sec:general} we prove the Theorem \ref{thm:general} using Steiner symmetrization and convexity.
In Section \ref{sec:maj} we prove Theorems \ref{thm:majn}, \ref{thm:maj1}, and \ref{thm:inverse-poly} --- the arguments in this section are mostly probabilistic.
Finally, in Section~\ref{sec:unproven} we discuss the results of some computer analysis, and pose two more open problems.

\subsection{Results obtained subsequent to this work}
Several new results have been proven about the cosmic coin problem
subsequent to this work. Perhaps the most interesting asymptotic
setting of the parameters is that of
Theorem~\ref{thm:inverse-poly}: $\eps$ fixed and $k \to \infty$,
with $n$ being unbounded.  It has since been shown that in this
setting the optimal success probability of the players is
precisely $\tilde{\Theta}(k^{-\nu})$, where $\nu = \nu(\eps) =
\frac{4\eps(1-\eps)}{(1-2\eps)^2}$, and $\tilde{\Theta}( \cdot )$
denotes asymptotics to within a subpolynomial ($k^{o(1)}$) factor.
The upper bound comes simply from a more careful analysis of the
success probability of the majority protocol.  Much more
interesting is the lower bound, which is proved using a
\emph{reverse} hypercontractive property of the Bonami-Beckner
operator proven in \cite{Bor82} 
(c.f.\ usual uses of the hypercontractive property in
\cite{KKL88} and subsequent works).

In addition, the same problem on different tree structures has
been studied. The problem in the present paper corresponds to a
star graph with $k$ leaves, with the initial random string at the
root distributed to the $k$ players along the edges.  The
corresponding problem on the line graph has subsequently been
studied, and it is shown that  
for $n$ unbounded and fixed $\eps$ the agreement probability decays
exponentially in $k$. This is different from the star where the decay
rate is polynomial in $k$.

See \cite{MORSS03} and also preliminary expositions from
\cite{O03}.

\section{Fourier methods} \label{sec:fourier}
In this section, we make a usual notational switch; the bits 0 and
1 will be denoted by $+1$ and $-1$, respectively. Note that $f : \fits^n \to \fits$
is balanced {\bf iff} $\E[f] = 0$.
Given the functions of the parties $f_1, \dots, f_k : \fits^n \to
\fits$, we view these as functions in the larger space $\fits^{n \times (k+1)}$ in the
natural way: $f_i(x, y^1, \dots, y^k) = f_i(y^i)$.
Our probability space gives rise to a natural inner product on
functions $f, g : \fits^{n \times (k+1)} \to {\bf R}$:
\begin{equation}
\langle f,g \rangle = \Ex_{x, y^1, \dots, y^k}[f(x,y^1,\dots,y^k)g(x,y^1,\dots,y^k)].
\end{equation}
\begin{lemma} \label{lem:fouriern}
\[
\sum_{i,j=1}^k \langle f_i, f_j \rangle = 2 \E[\#(i,j) : f_i(y^i) = f_j(y^j)] - k^2.
\]
\end{lemma}

\begin{proof}
Since the $f_i$ are $\pm 1$ valued functions,
\begin{eqnarray*}
\sum_{i,j=1}^k \langle f_i, f_j \rangle
&=& \sum_{i,j=1}^k \left( \P[f_i(y^i)  = f_j(y^j)] - \P[f_i(y^i)  \neq f_j(y^j)] \right) \\
&=& \sum_{i,j=1}^k \left( 2 \P[f_i(y^i)  = f_j(y^j)] - 1 \right) = 2 \E[\#(i,j) : f_i(y^i) = f_j(y^j)] - k^2.
\end{eqnarray*}
\end{proof}

Now, in order to maximize the quantity in (\ref{eq:nagreements}), we analyze the scalar products $\langle f_i,f_j \rangle$.
In order to analyze scalar products, it is useful to work with the Fourier basis. We refer the reader to
\cite{KKL88,BKS99,MO02} for background.
For a set $S \subset [n]$, we let $U_{i,S} : \fits^{n \times (k+1)} \to \fits$ be defined by:
\[
U_{i,S}(x,y^1,\ldots,y^k) = U_{i,S}(y^i) = \prod_{j \in S} y^i_j.
\]
Since the $k$ $y_j$'s are independent, it follows that if $S \neq
S'$ then for all $i,i'$,
\begin{equation} \label{eq:Uorth}
\langle U_{i,S}, U_{i',S'} \rangle = 0.
\end{equation}
Moreover, if $i \neq i'$ then
\begin{equation} \label{eq:Unoisorth}
\langle U_{i,S}, U_{i',S} \rangle = \E[ \prod_{j \in S} y^{i}_j y^{i'}_j] = \prod_{j \in S} \E[y^{i}_j y^{i'}_j]
 = (1 - 2 \eps)^{2 |S|},
\end{equation}
and the functions $U_{i,S}$ all have norm $1$, so $\langle U_{i,S}, U_{i,S} \rangle = 1$.

\begin{lemma} \label{lem:k=2}
Let $i \neq j$.
Then
\[
\max_{f_i, f_j \in \B_n} \langle f_i(y^i), f_j(y^j) \rangle = (1 - 2 \eps)^2,
\]
and the maximum is obtained when $f_i = f_j = f$ and $f(x) = \pm x_r$ for some $1 \leq r \leq n$.
\end{lemma}

\begin{proof}
Express $f_i$ and $f_j$ in terms of their Fourier expansion, $f_i
= \sum_{S \subseteq [n]} \hat{f}_i(S) U_{i,S}$ and similarly for $f_j$.
Since, both $f_i$ and $f_j$ are balanced, and $\E[U_{i,S}] = \E[U_{j,S}] = 0$ for nonempty $S$,
it follows that $\hat{f}_i(\emptyset) = \hat{f}_j(\emptyset) = 0$.
Now by (\ref{eq:Uorth}) and (\ref{eq:Unoisorth}) it follows that
\begin{equation} \label{eqn:fourier}
\langle f_i,f_j \rangle =
\sum_{\emptyset \neq S \subseteq [n]} \hat{f}_i(S) \hat{f}_j(S) (1-2\epsilon)^{2|S|}.
\end{equation}
Hence we have:
\begin{eqnarray*}
\langle f_i, f_j \rangle & = & \sum_{\emptyset \neq S \subseteq
[n]} (1-2\epsilon)^{|S|}\hat{f}_i(S)
(1-2\epsilon)^{|S|}\hat{f}_j(S) \\
& \leq & \sqrt{\sum_{\emptyset \neq S \subseteq [n]}
(1-2\epsilon)^{2|S|} \hat{f}_i(S)^2} \sqrt{\sum_{\emptyset \neq
S \subseteq [n]} (1-2\epsilon)^{2|S|} \hat{f}_j(S)^2}
\qquad\text{(Cauchy-Schwarz)} \\
& \leq & \sqrt{\sum_{\emptyset \neq S \subseteq [n]}
(1-2\epsilon)^2 \hat{f}_i(S)^2} \sqrt{\sum_{\emptyset \neq S
\subseteq [n]} (1-2\epsilon)^2 \hat{f}_j(S)^2} \\
& = & (1-2\epsilon)^2,
\end{eqnarray*}
as $\sum \hat{f}_i(S)^2 = \sum \hat{f}_j(S)^2 = 1$.
The second inequality is tight only if $f_i$ and $f_j$ have Fourier degree 1.
Note that if $f(x) = \sum_{|S|=1} \hat{f}(S) u_S(x)$ is a function which is $\pm 1$ valued, then for all $S$ of size $1$,
it holds that $2 \hat{f}(S) = f(x) - f(x \oplus e_S) \in \{-2,0,2\}$. It follows that $f(x) = \pm x_r$ for some $r$.

In this case, the first inequality is tight only if $f_i$ and $f_j$ are the same one-bit function.  Hence, as claimed, $f_i = f_j = f$ where $f(x) = \pm x_r$ constitutes the only maximizing solution.
\end{proof}

We can now prove Theorem~\ref{thm:fourier}.

\begin{proof} [of Theorem \ref{thm:fourier}]
By Lemma \ref{lem:fouriern} it follows that maximizing $\E[\#(i,j) : f_i(y^i) = f_j(y^j)]$ is the same as maximizing,
\[
\sum_{i,j=1}^k \langle f_i, f_j \rangle = k + \sum_{i \neq j} \langle f_i, f_j \rangle.
\]
By Lemma \ref{lem:k=2}, the above sum is maximized when $f_1 = \ldots = f_k = f$, and $f(x) = x_1$ up to $\pi_S$.  We thus obtain the first assertion of the theorem.

For the second assertion, note that when $k=2$,
\[
\E[\#(i,j) :  f_i(y^i) = f_j(y^j)] = 2 + 2 \P[f_1(y^1) = f_2(y^2)],
\]
while when $k=3$,
\[
\E[\#(i,j) : f_i(y^i) = f_j(y^j)] = 5 + 4 \P[f_1(y^1) = f_2(y^2) = f_3(y^3)],
\]
so the second assertion follows.
\end{proof}

\section{Convexity and symmetrization} \label{sec:general}
We now show that to maximize $\p(f_1, \dots, f_k ; \epsilon)$, it suffices to look at restricted sets of functions.
The methods in the section are related to convexity in general and the Steiner symmetrization in particular, see e.g.\@ \cite{T93} for background.

We begin with a definition and a simple Fourier Lemma (The
Bonami-Beckner operator $T_{\eps}(f)$ was first defined in \cite{Be75,Bo70},
see also \cite{KKL88,BKS99})
\begin{definition} \label{def:beckner}
For $f : \bits^n \to \R$, given by $f(x) = \sum_{S \subset [n]} \hat{f}(S) u_S(x)$,
let
\[
T_{\eps}(f)(x) = \sum_{S \subset [n]} \hat{f}(S) (1 - 2 \eps)^{|S|} u_S(x).
\]
\end{definition}
\begin{lemma}\label{lem:invert}
Given $\eps$ and $f : \bits^n \to \bits$,
$T_\eps(f)(x)$ equal the probability that a particular party using $f$ outputs 1, given that the source string is $x$.
If $f$ and $g$ are different boolean functions on $\bits^n$, then for every $0 < \eps < 1/2$,
there exists some $x \in \bits^n$ for which $T_\eps(f)(x) \neq T_\eps(g)(x)$.
\end{lemma}
\begin{proof}
Note that given a source $x$ and an $\eps$ corrupted version of
$x$, $y$, the expected value of $u_S(y)$ is given by $(1 -
2\eps)^{|S|} u_S(x)$ (an easy calculation; see e.g. \cite{BKS99}
for a full proof). Therefore, by linearity of expectation, it
follows that for all $f : \bits^n \to \R$, $T_{\eps}(f)(x)$ is the
expected value of $f(y)$, where $y$ is a corrupted version of $x$.
In particular, if $f : \bits^n \to \bits$, then $T_{\eps}(f)(x) =
\E[f(y)] = \P[f(y) = 1]$ and we proved the first assertion of the
lemma.

For the second assertion note that for $0 < \eps < 1/2$, $T_{\eps}$ is by definition a reversible linear transformation
on the space of all function from $\bits^n \to \R$.
\end{proof}

Now, using convexity, we prove that all parties should use the same function:
\begin{proposition}
\label{prop:same} Fix $k$, $n$, and $0 < \eps < 1/2$.  Let $\C$
be any class of boolean functions on $n$ bits.  Subject to the
restriction that $f_1, \dots, f_k \in \C$, every protocol which maximizes
$\p(f_1, \dots, f_k ; \epsilon)$ has $f_1 = \cdots = f_k$.
\end{proposition}
\begin{proof}
Let $\C = \{f_1, f_2, \dots, f_M\}$, and assume $M > 1$ -- else the proposition is trivial.
Suppose that among the $k$ parties, exactly $t_j$ use the function $f_j$.  Then clearly,
\begin{equation}
t_j \geq 0, \qquad \sum_{j=1}^M t_j = k, \qquad t_j \in \Z.  \label{eq:same}
\end{equation}
The probability that all parties agree is:
\begin{equation}
\p = \sum_{x \in \bits^n} 2^{-n} \left(\prod_{j=1}^M
\left(T_\eps(f_j)(x)\right)^{t_j} + \prod_{j=1}^M \left(1-T_\eps(f_j)(x)\right)^{t_j} \right).
\end{equation}
Note that each $T_\eps(f_j)(x) \in (0,1)$ and that for any $c \in (0,1)$, the
function $g(t) = c^t$, is {\em log-convex} (since $\log c^t = t \log c$ is linear).
Therefore the function $g_1 \cdots g_M : \R^M \to \R$ given by $(t_1, \dots,
t_M) \mapsto \prod_{j=1}^M g_j(t_j)$ is a log convex function, and therefore
a convex function.

Since the sum of convex functions is also convex, $\p$ is a convex function of the $t_j$'s.
We wish to maximize $\p$ subject to the restrictions (\ref{eq:same}).

If we relax the assumption $t_j \in \Z$ to $t_j \in \R$,
we are simply maximizing a convex function over a convex bounded polytope.
The vertices of the polytope are simply the points of the form $(0, \dots, 0, k, 0, \dots, 0)$.
The maximum must occur at a vertex,
and so it follows that there is at least one maximizing protocol in which all players use the same function.

It remains to show that $\p$ doesn't obtain the maximum at any point which is not a vertex of the polytope.
Note that by convexity, if $\p$ has a maximum which is not a vertex of the polytope, then there exists an interval
$I = \{ t v_1 + (1 - t) v_2 : t \in [0,1]\}$, where $v_1$ and $v_2$ are vertices of the polytope, such that $f$
is a constant function on $I$. Therefore if we could show that $f$ is strictly convex on $I$ (as a function of $t$),
then it will follow that the maximum is obtained only at vertices of the polytope.

note that when restricted to the edge $I$ joining, e.g., $v_1 = (k, 0, \dots, 0)$ and $v_2 = (0, k, 0, \dots, 0)$,
$\p$ is given by:
\[\p = \sum_{x \in \bits^n} 2^{-n} \left(\left(T_\eps(f_1)(x)\right)^{\lambda k} \left(T_\eps(f_2)(x)\right)^{(1-\lambda) k} + \left(1-T_\eps(f_1)(x)\right)^{\lambda k} \left(1-T_\eps(f_2)(x)\right)^{(1-\lambda) k}\right).
\]
By Lemma~\ref{lem:invert}, we can find an $x_0$ such that $T_\eps(f_1)(x_0)$ and $T_\eps(f_2)(x_0)$ differ.
Therefore the function
\[
\left(T_\eps(f_1)(x_0)\right)^{\lambda k} \left(T_\eps(f_2)(x_0)\right)^{(1-\lambda) k} =
\left(T_{\eps}(f_2)(x_0)\right)^k \left(\frac{T_\eps(f_1)(x_0)}{T_\eps(f_2)(x_0)} \right)^{\lambda k},
\]
is strictly convex, and $\p$ is strictly convex on $I$ as needed.
\end{proof}

Next we use Steiner symmetrization principle in order to obtain more information on functions which
optimize $\p(f_1, \dots, f_k ;\epsilon)$.
Recall that for $x, y \in \bits^n$, we write $x \preceq y$ if for all $i \in [n]$ it holds that $x_i \leq y_i$, and
we say that $f$ is monotone if $f(x) \leq f(y)$, whenever $x \preceq y$.
Similarly for $S \subset [n]$, we may write $x \preceq_S y$, if $x_i \leq y_i$ for $i \in S$, and
$y_i \leq x_i$ for $i \notin S$. We call a function $f$ which is monotone with respect to $\preceq_S$, $S$-monotone.
\begin{proposition}
\label{prop:monotone} Let $\C$ stand for either $\B_n$ or $A_n$.
For any $k$, $n$, $\epsilon$, if $f$ is restricted to be in $\C$,
and the maximum of $\p_k(g; \epsilon)$ is obtained at $g=f$, then
$f$ is $S$-monotone for some set $S$. Moreover, there exists $f
\in \C$ which maximizes $\p_k(g; \epsilon)$ and is monotone.
\end{proposition}
\begin{proof}
Let $f \in \C$ be any function which maximizes $\p_k(f; \epsilon)$ among functions in $\C$.  Let $f'$ be obtained from $f$ by ``shifting'' up in the first coordinate: Given $x \in \bits^{n-1}$,
\begin{itemize}
\item if $f(0x) = f(1x)$, then set $f'(0x) = f'(1x) = f(0x) = f(1x)$;
\item if $f(0x) \neq f(1x)$ then set $f'(0x) = 0$, $f'(1x) = 1$.
\end{itemize}
It is easy to see that in the case $\C = \B_n$, $f'$ remains in $\C$; a little thought reveals that this is again true in the case $\C = \A_n$.

For $y \in \bits^n$, let $\tilde{y} \in \bits^{n-1}$ be the last $n-1$ bits of $y$.
We claim that $\p_k(f'; \epsilon) \geq \p(f; \epsilon)$.  To show this, it suffices to show that for all $z^1,\ldots,z^k \in \bits^{n-1}$, \begin{equation} \label{eq:compshift1}
\P[f'(y^1) = \ldots = f'(y^k) \; | \; \tilde{y}^1 = z^1, \ldots, \tilde{y}^k = z^k] \geq
\P[f(y^1) = \ldots = f(y^k)  \; | \;  \tilde{y}^1 = z^1, \ldots, \tilde{y}^k = z^k].
\end{equation}
So suppose each $y^i$'s last $n-1$ bits are fixed to be $z^i$. Given $z^i$, $f(y^i)$ is a function from $\bits$ to $\bits$, and is therefore either the constant function $0$,
the constant function $1$, the identity function $id$, or the function $x \to \bar{x}$, which we denote by $\overline{id}$.

If $f(y^i)$ is already determined by $z^i$, then so is $f'(y^i)$ and the determined value is the same.  Otherwise, $f(y^i)$ is a function of the one remaining unknown bit,
$y^i_1$, and is either the function $id$ or $\overline{id}$. In {\em either} case, $f'(y^i)$ is the identity function on $y^i_1$.

Assume that given $(z^1,\ldots,z^k)$, there are $a+b$ undetermined functions $f(y^i_1)$, with $a$ of them $id$, and $b$ of them $\overline{id}$.
The probability that all of these functions agree on $0$ (or $1$) is
\[
q = \frac{1}{2} \left( (1-\epsilon)^a \epsilon^b + \epsilon^a(1-\epsilon)^b \right),
\]
and the probability that all of the undetermined $f'$'s agree on 0 (or $1$) is
\[
q' = \frac{1}{2} \left( (1-\epsilon)^{a+b} + \epsilon^{a+b} \right).
\]

There are three cases to consider:
\begin{itemize}
\item
If some of the determined functions are determined to be $0$ and some to be $1$, then both terms in (\ref{eq:compshift1}) are zero.
\item
If all of the determined functions are determined to be $0$ ($1$), then the left side of (\ref{eq:compshift1}) is $q'$ and the right side of (\ref{eq:compshift1}) is $q$.
\item
If there are no determined functions, then the left side of (\ref{eq:compshift1}) is $2 q'$ and the right side of (\ref{eq:compshift1}) is $2 q$.
\end{itemize}

Therefore the claim will follow once we show that $q' \geq q$.
\begin{eqnarray}
\frac{1}{2} \left((1-\epsilon)^{a+b} + \epsilon^{a+b} \right) & \geq &
\frac{1}{2} \left((1-\epsilon)^a \epsilon^b + \epsilon^a(1-\epsilon)^b \right)\\
\Leftrightarrow \qquad 1 + \left(\frac{\epsilon}{1-\epsilon}\right)^{a+b} &
\geq & \left(\frac{\epsilon}{1-\epsilon}\right)^b +
\left(\frac{\epsilon}{1-\epsilon}\right)^a, \label{eqn:strict}
\end{eqnarray}
which follows by the convexity of the function $t \to
\left(\frac{\epsilon}{1-\epsilon}\right)^t$.

Thus we've established $\p_k(f'; \epsilon) \geq \p(f; \epsilon)$.
We further claim that this inequality is strict unless $f$ was already monotone or anti-monotone on the first coordinate.
If $f$ is neither monotone nor anti-monotone on the first coordinate,
then there exist $z^1$ and $z^2$ such that $f$, when the last $n-1$ coordinates are restricted to $z^1$,
becomes $id$, and when the last $n-1$ coordinates are restricted to $z^2$, becomes $\overline{id}$.
Picking $z^3, \dots, z^k$ so that all the other restricted functions are either $id$ or $\overline{id}$,
we obtain $a, b \geq 1$, so (\ref{eqn:strict}) is strict inequality and therefore $q' > q$.

Repeating the above argument for all other coordinates, it follows that any maximizing function $f$ must be
$S$ monotone and that there exists a maximizing function which is monotone.
\end{proof}

Recall that for $x,y \in \bits^n$, we write $x \preceq_L y$ if
$\sum_{i=1}^m x_i \leq \sum_{i=1}^m y_i$ for every $m=1\ldots n$, and that
we call $f : \bits^n \to \bits$ left-monotone, if $f(x) \leq f(y)$ whenever $x \preceq_L y$.
\begin{proposition} \label{prop:left-monotone}
Let $\C$ stand for either $\B_n$ or $A_n$.
Let $k$, $n$, and $0 < \epsilon < 1/2$.
Suppose that $f$ maximizes $\p_k(f; \epsilon)$ in $\C$, then up to $\pi_S$, $f$ is left-monotone.
\end{proposition}

\begin{proof}
The proof is similar to the proof of Proposition~\ref{prop:monotone}, so we will be more brief.
By Proposition~\ref{prop:monotone}, we may assume that $f$ is monotone.

Now apply a new sort of shift to $f$.  Suppose we fix all but two input bits
to $f$.  Since $f$ is monotone, there are only $6$ possibilities for
what the restricted function is; its support may be $\emptyset$,
$\{11\}$, $\{11,10\}$, $\{11,01\}$, $\{11,10,01\}$ or $\{11,10,01,00\}$.
Define $f'$ to be the same function in all cases except when the
support is $\{11,01\}$; in this case, switch it to $\{11,10\}$.
This rule preserves balance and asymmetry.

We want to show that $\p_k(f' ;  \epsilon) \geq \p_k(f; \eps)$.
As before, we condition on all but two bits of each
of $y_1, \dots, y_k$, and show that $f'$ is better.  Say
that under this conditioning, $a$ of the $f(y^i)$'s restrict to the
function with support $\{11,10\}$, and $b$ of the $f(y^i)$'s restrict to the function with support $\{11,01\}$.  Since all other possible restricted functions have the same value for $01$ as they do for $10$, it suffices again to compare the probability
with which the $a+b$ functions agree on $1$ with the probability that
the corresponding shifted functions agree on $1$. Further, by symmetry, we need
only consider the cases when the two source bits from $x$ are
different (otherwise $f$ and $f'$ do equally well).

So considering the two cases --- the source bits are $10$ or the source bits are $01$ --- we get that the contribution from the $f$-restricted functions will be $(1/2)((1-\epsilon)^a \epsilon^b) +
\epsilon^b(1-\epsilon)^a)$, and the contribution from their shifted versions
will be $(1/2)((1-\epsilon)^{a+b} +
\epsilon^{a+b})$.  As we saw in Proposition~\ref{prop:monotone}, this
latter quantity is always at least the former quantity.  Hence the
shift can only improve the probability of agreement.

Hence we indeed have $\p_k(f'; \eps) \geq \p_k(f; \eps)$.
If we repeatedly apply this shift to all pairs of coordinates, we end up with a left-monotone function.

Note that if none of the shifting operations strictly increased the probability of agreement for $f$,
then for every pair of coordinates $(i,j)$ which were shifted, either all the balanced restrictions of $f$ to coordinates
$(i,j)$ have support $\{11,10\}$, or all the balanced restrictions have support $\{11,01\}$.
In either case, all the shifting did is to replace the functions $f$ by a function $f \circ \pi_{\emptyset}$,
where $\pi_{\emptyset}$ is the transposition of coordinates $(i,j)$. It thus follows that the original function
was left monotone up to some $\pi_{\emptyset}$, as needed.
\end{proof}

\begin{proof} [of Theorem \ref{thm:general}]
The proof follows from Propositions \ref{prop:same}, \ref{prop:monotone}, and \ref{prop:left-monotone}.
\end{proof}


\section{Majorities} \label{sec:maj}
In this section we study majority functions and show that these function are optimal for some limiting values of $k$ and $\eps$.

\subsection{Fixed $\eps$, $n$; $k \to \infty$}
We start by proving Theorem \ref{thm:majn}.
Given a function $f : \bits^n \to \bits$, let $p_1(f,x,\epsilon)$
denote the probability that $f(y) = 1$, given that $y$ is an
$\epsilon$-corrupted version of the string $x \in \bits^n$.
Let $p_0$ be defined similarly for the probability that $f(y) = 0$.
\begin{proposition} \label{prop:1-maximizes}
Fix $\eps$, and let $f$ be monotone.  Then as a
function of $x$, $p_1(f,x,\epsilon)$ is maximized at $x = \vec{1} = 1\cdots 1$,
and $p_0(f,x,\epsilon)$ is maximized at $x = \vec{0} = 0 \cdots 0$.
\end{proposition}
\begin{proof}
We prove the claim for $p_1$, the proof for $p_0$ being the same.
Note that flipping each bit of a string with probability
$\epsilon$ is the same as {\em updating} each bit with probability
$2\epsilon$, where an update consists of replacing the bit with a
random choice from $\bits$.

Let $x \in \bits^n$ be any sequence.  Let $x'$ be an $\eps$ corrupted version of $x$, and $\vec{1}'$ be an $\eps$ corrupted version of $\vec{1}$.
We claim that we can couple the random variable $x'$ and $\vec{1}'$
in such a way that $x' \preceq \vec{1}'$.

The coupling is achieved in the following simple way: update the same bits
of $x$ and $\vec{1}$ with the same values.
Clearly, we have $x' \preceq \vec{1}'$.  Hence by
monotonicity if $f(x') = 1$, then $f(\vec{1}') = 1$.  The result
follows, as
\[
p_1(f,x,\epsilon) = \P[f(x') = 1] \leq \P[f(\vec{1}') = 1] = p_1(f,\vec{1},\eps).
\]
\end{proof}

\begin{proposition}
\label{prop:maj-best-1} For fixed $\eps$, $p_1(f,\vec{1},\epsilon)$ and $p_0(f,\vec{0},\epsilon)$
are maximized among $f \in \B_n$ by any function which is 0 on all strings with fewer than $n/2$ 1's.
In particular, if $n$ is odd, $f = \maj{n}$ is the unique maximizing function.
\end{proposition}
\begin{proof}
We prove the assertion about $p_1$.
\[
p_1(f, \vec{1}, \epsilon) = \sum_{x \in f^{-1}(1)}
(1-\epsilon)^{n- \Delta(x, \vec{1})} \epsilon^{\Delta(x, \vec{1})},
\]
where $\Delta$ denotes Hamming distance.  The quantity being summed is a strictly decreasing function of $\Delta(x, \vec{1})$.  The result follows.
\end{proof}

\begin{proof} [of Theorem \ref{thm:majn}]
We prove the theorem for $n$ odd. The proof for $n$ even is essentially the same.
By Theorem \ref{thm:general}, we may assume without loss of generality
that all parties use the same monotone function $f \in \B_n$.  Now:
\begin{eqnarray} \label{eq:majexp}
\p(\maj{n}, k, \epsilon) & = & 2^{-n}\sum_{x \in \bits^n}
\left( p_1(\maj{n}, x, \epsilon)^k + (p_0(\maj{n}, x, \epsilon))^k \right) \\ \nonumber
& \geq & 2^{-n} \left( p_1(\maj{n}, \vec{1}, \epsilon)^k + p_0(\maj{n}, \vec{0}, \epsilon)^k \right).
\end{eqnarray}
By Proposition~\ref{prop:1-maximizes}, if $f$ is monotone, then
\begin{eqnarray} \label{eq:montexp}
\p(f, k, \epsilon) & = & 2^{-n}\sum_{x \in \bits^n}
\left(p_1(f, x, \epsilon)^k + (p_0(f, x, \epsilon))^k \right) \\ \nonumber
& \leq & 2^{-n}\sum_{x \in \bits^n} \left( p_1(f, \vec{1}, \epsilon)^k + p_0(f, \vec{0}, \epsilon)^k \right)
= p_1(f, \vec{1}, \epsilon)^k + p_0(f, \vec{0}, \epsilon)^k.
\end{eqnarray}
By Proposition~\ref{prop:maj-best-1}, if $f \in \B_n$ is monotone,
and $f \neq \maj{n}$, then $p_1(f, \vec{1}, \epsilon) < p_1(\maj{n}, \vec{1}, \epsilon)$ and
$p_0(f, \vec{0}, \epsilon) < p_0(\maj{n}, \vec{0}, \epsilon)$.
Hence for sufficiently large $k$, we will have
\begin{equation} \label{eq:majmont}
\begin{array}{ll}
2^{-n} p_1(\maj{n}, \vec{1}, \epsilon)^k > p_1(f, \vec{1},\epsilon)^k, &
2^{-n} p_0(\maj{n}, \vec{0}, \epsilon)^k > p_0(f, \vec{0},\epsilon)^k.
\end{array}
\end{equation}
Combining (\ref{eq:majexp}), (\ref{eq:montexp}) and (\ref{eq:majmont})
we obtain that $\p_k (\maj{n} ; \epsilon) \geq \p_k(f ; \epsilon)$ for all monotone $f \in \B_n$ as needed.

Bound (\ref{eq:as_maj}) follows from (\ref{eq:majexp}), (\ref{eq:montexp}), and (\ref{eq:majmont}) once we note that
\[
p_0(\maj{n}, \vec{0}, \epsilon) =
p_1(\maj{n}, \vec{1}, \epsilon) = 1 - \P[Bin(n,\eps) > n/2].
\]
\end{proof}

\subsection{Fixed $k$, $n$; $\eps \to 0$ or $1/2$}
\begin{proposition} \label{prop:eps0}
For all $k$ and $n$, there exists $\eps'(k,n) > 0$, such that for all $0 < \eps < \eps'(k,n)$, the unique best protocol up to $\pi_S$ for maximizing $\p(f_1,\ldots,f_k ; \eps)$ is given by $f_1 = \ldots = f_k = f$, where $f(x) = x_1$.
\end{proposition}

\begin{proof}
From Proposition~\ref{prop:same}, it follows that the maximum can only be obtained if
$f_1 = \ldots = f_k = f$.
Note that the probability that there is more than one corrupted bit is $O(\eps^2)$ (the constant in the $O(\cdot)$ {\em does} depend on $k$ and $n$).
Suppose that only the $i$th bit for party $j$ was corrupted. Then all the parties will agree if and only if
$f(x) = f(x \oplus e_i)$, where $x \oplus e_i$ is the vector $x$ with the $i$th bit flipped.
We therefore obtain,
\begin{equation} \label{eq:ap1}
\p_k(f; \eps) = (1 - \eps)^{k n} + k \eps (1 - \eps)^{n k - 1} \sum_{i=1}^n \P_x[f(x) = f(x \oplus e_i)] + O(\eps^2).
\end{equation}
Writing $A$ for the set $\{x : f(x) = 1\}$, and $\partial_E(A)$ for the {\em edge-boundary} of the set $A$,
\[
\partial_E(A) = \cup_{i=1}^n \{(x, x \oplus e_i) : x \in A, x \oplus e_i \notin A\},
\]
we see that
\[
\sum_{i=1}^n \P_x[f(x) = f(x \oplus e_i))] = n - 2^{-n+1} \partial_E(A).
\]
So for small $\eps$ (compared to $k$ and $n$), in order to maximize $\p_k(f ; \eps)$, we should minimize $\partial_E(A)$ for sets such that
$|A| = 2^{n-1}$. By the isoperimetric inequality for the cube, the sets $A$ which minimize $\partial_E(A)$ among all
sets of size $2^{n-1}$ are exactly the sets $A = \{ x : x_i = 0\}$, or $A = \{x : x_i = 1\}$. Thus $f$ must be $f(x) = x_1$ up to $\pi_S$, as claimed.
\end{proof}

\begin{proposition} \label{prop:eps1/2}
For all $k$ and $n$, there exists $\eps'(k,n) < 1/2$, such that for all $eps'(k,n) < \eps < 1/2$, the unique best protocol up to $\pi_S$ for maximizing $\p(f_1,\ldots,f_k ; \eps)$ is given by $f_1 = \ldots = f_k = f$, where $f(x) = x_1$.
\end{proposition}

\begin{proof}
Again, Proposition~\ref{prop:same} implies that we need only consider the case $f_1 = \ldots = f_k = f$.
In this proof it will be helpful again to assume that the bit values are $\pm 1$, so we want to show that
the maximizing functions are $f(x) = x_i$, or $f(x) = -x_i$.

It will be useful to work with the ``updating representation''.
Let $X(i,j)$ for $1 \leq i \leq k$ and $1 \leq j \leq n$ be a sequence of i.i.d $\{0,1\}$ variables s.t.\@ $\P[X(i,j) = 1] = \delta = 1 - 2 \eps$. Note that we may produce the $y^i$'s from $x$ in the following manner.
If $X(i,j) = 1$, then $y^i_j = x_j$, otherwise  $y^i_j$ is chosen uniformly at random from $\fits$ independently from everything else.

Note that if all the $X(i,j)$ are 0, the inputs to the functions are independent, so for all balanced
$f$'s,
\[
\P[f(y^1) = \ldots = f(y^k)\; | \; \sum_{i,j} X_{i,j} = 0] = 2^{-k+1}.
\]
Similarly for all balanced $f$,
\[
\P[f(y^1) = \ldots = f(y^k) \; | \; \sum_{i,j} X_{i,j} = 1] = 2^{-k+1},
\]
and for all $i,i'$ and $j \neq j'$,
\[
\P[f(y^1) = \ldots = f(y^k) \; | \; X_{i,j} = X_{i',j'} = 1, \sum_{s,t} X_{s,t} = 2] = 2^{-k+1}.
\]
Moreover,
\[
\P\left[\sum_{i,j} X_{i,j} > 2\right] = O(\delta^3)
\]
(the constant in the $O(\cdot)$ depending on $k$ and $n$).

We therefore conclude that
\begin{equation} \label{eq:d3}
\p_k(f ; \eps) = c_k + \delta^2 (1 - \delta)^{n k - 2} \sum_{i \neq i'} \sum_{j=1}^n
\P[f(y^1) = \ldots = f(y^k)\;| \; X_{i,j} = X_{i',j} = 1] + O(\delta^3),
\end{equation}
where $c_k$ is independent of $f$.
Writing $z$ for a uniformly chosen element of $\fits^n$, and $z'$ for an element which is chosen by picking $1 \leq i \leq n$
uniformly at random, and then choosing $z' \in \bits$ uniformly among all $z'$ s.t. $z'_i = z_i$, we obtain,
\begin{equation} \label{eq:zz'}
\sum_{i \neq i'} \sum_{j=1}^n
\P[f(y^1) = \ldots = f(y^k)\; | \; X_{i,j} = X_{i',j} = 1] = n k(k-1) \P[f(z) = f(z')].
\end{equation}
Therefore, if we could show that $\P[f(z) = f(z')]$ is maximized among all balanced functions $f$ when $f(x) = x_1$ up to $\pi_S$, then the proof will follow from (\ref{eq:d3}).

In order to prove this claim, we note that
\[
\P[f(z) = f(z')] = (1 + \E[f(z) f(z')])/2.
\]
It therefore suffices to maximize $\E[f(z) f(z')]$ over all balanced functions $f$, i.e., functions with $\hat{f}(\emptyset) = \E[f] = 0$.
We pass to the Fourier representation as in the proof of Theorem \ref{thm:fourier}\@. It is easy to see that if
$u_S(x) = \prod_{i \in S} x_i$, then
\[
\E[u_S(z') | z] = \left\{ \begin{array}{ll} u_S(z) & \mbox{ if } S = \emptyset, \\
                                            u_S(z)/n & \mbox{ if } |S| = 1, \\
                                            0          & \mbox{otherwise.} \end{array}
                                            \right.
\]
So,
\[
\E[f(z) f(z')] = \frac{1}{n}\sum_{|S| = 1} \hat{f}^2(S) \leq 1/n,
\]
and equality is achieved iff $f(x) = x_i$, or $f(x) = -x_i$ as needed (see proof of Lemma \ref{lem:k=2}).
\end{proof}

\begin{proof} [of Theorem \ref{thm:maj1}] Follows immediately from Propositions~\ref{prop:eps0} and \ref{prop:eps1/2}.
\end{proof}

\subsection{Fixed $\eps$; $k \to \infty$ with $n$ unbounded}
Finally, we prove Theorem~\ref{thm:inverse-poly}.

\begin{proof} [of Theorem~\ref{thm:inverse-poly}]
Fix $\eps$ and $k$.  We consider $\p_k(\maj{n}; \eps)$ as a function of $n$, as $n \to \infty$ through the odd numbers.  Our proof will go by showing that there is at least a $\Omega(1/k^{2.01/(1-2\eps)^2})$ chance that the source string $x$ has significantly more 1's than 0's.  Then we show that in this case, the probability any particular party says 1 is at least $1 - 1/k$, and hence the probability that all parties say 1 is at least a constant.

Let $X$ be the random variable given by $(\text{\# 1's in } x) - (\text{\# 0's in } x)$.  By the Central Limit Theorem, as $n \to \infty$, the distribution of $X$ approaches a normal distribution with mean 0 and variance $n$.  Let $c = \frac{2}{1 - 2\eps}$.  The probability that an $N(0,n)$ normal variable
exceeds $c \sqrt{\log k} \sqrt{n}$ is:
\[
1 - \Phi(c\sqrt{\log k}) \geq \frac{1}{2}\frac{1}{c\sqrt{\log k}}
\frac{1}{\sqrt{2\pi}}\exp\left(-c^2 \log k / 2\right) \geq \Omega(1/k^{2.01/(1-2\eps)^2}).
\]
Here $\Phi$ denotes the cumulative distribution function of a
standard normal variable, and the first inequality follows from
the fact that $1 - \Phi(x) \geq (1/x - 1/x^3)\phi(x)$ (see
\cite{F68}), where $\phi(x)$ is the density function of a standard
normal variable .

Given that this happens, pessimistically assume that $x$ contains
just $c \sqrt{\log k} \sqrt{n}$ more 1's than 0's; i.e., $x$
contains exactly $n/2 + (c/2) \sqrt{\log k}\sqrt{n}$ 1's.  We now
show that the probability that a particular party using $\maj{n}$
outputs 1 given $x$ is at least $1 - 1/k$.

Consider the $\epsilon$-corrupted version of $x$ the party sees;
call it $y$.  The number of 1's in $y$ is distributed as the sum of $n$ Bernoulli trials, $n/2 + (c/2)\sqrt{\log k}\sqrt{n}$ of which have success probability $1-\eps$, and $n/2 - (c/2)\sqrt{\log k}\sqrt{n}$ of which have success probability $\eps$.  We can use a single Chernoff bound to upper-bound the probability of getting fewer than $n/2$ 1's in $y$.  The expected number of 1's is $n/2 + (1-2\eps)c\sqrt{\log k}\sqrt n = n/2 + 2\sqrt{\log k}\sqrt n$.  Since $n/2 = (1 - \delta)(n/2 + 2\sqrt{\log k}\sqrt n)$ when $\delta = 2\sqrt{\log k}\sqrt n/(n/2 + 2\sqrt{\log k}\sqrt n) > 2\sqrt{\log k}/\sqrt{n}$, Chernoff tells us that the probability that $y$ has fewer than $n/2$ 1's is at most $\exp(-\frac{4\log k\;\;(n/2)}{2n}) = 1/k$.

Thus as claimed, given a source string with at least $c\sqrt{\log k}\sqrt{n}$ more 1's than 0's, the probability a particular party outputs 1 is at least $1-1/k$.  Hence the probability that all parties output 1 is at least $(1-1/k)^k = \Omega(1)$.

Hence by taking $n$ sufficiently large, we can make $\p_k(\maj{n}; \epsilon) \geq \Omega(1/k^{2.01/(1-2\eps)^2})$.  By applying the Berry-Ess\'een bounds on the rate of convergence in the Central Limit Theorem (see \cite{F68}), one can show that it suffices for $n$ to be $O(k^{4.01/(1-2\eps)^2})$.
\end{proof}

In the limit as $n \to \infty$, all distributions involved in the calculation of $\p_k(\maj{n}; \eps)$ become normal, and it is possible to get some more or less closed forms for the limit:
\begin{proposition}\label{prop:integral}
\begin{eqnarray}
\lim_{\substack{n \to \infty \\ n\text{ odd}}} \p_k(\maj{n}; \eps) &=& 2\,\E\left[\Phi\left(\frac{X}{\sqrt{c(\eps)}}\right)^k\right] \label{eqn:ex}\\
&=& \frac{2\sqrt{c(\eps)}}{(2\pi)^{\frac12(c(\eps) - 1)}} \int_0^1 x^k I(x)^{c(\eps)-1}\,\mathrm{d}x,\label{eqn:int}
\end{eqnarray}
where $c(\eps) := \frac{4\eps(1-\eps)}{(1-2\eps)^2} \in (0,\infty)$, $X$ is a standard normal random variable, and $I(x) := \phi(\Phi^{-1}(x))$.
\end{proposition}
(We thank Nati Srebro for his help in calculating (\ref{eqn:ex}).)  The second formula (\ref{eqn:int}) can be used to get tighter bounds than in Theorem~\ref{thm:inverse-poly}.  For example, from (\ref{eqn:int}) we get that $\lim_{n \to \infty} \p_k(\maj{n}; 1/2 - \sqrt{2}/4) = 2/(k+1)$.

\section{Computer-assisted results and open problems} \label{sec:unproven}
\label{sec:computer}
The problem well avails itself to analysis by computer.  In particular, given any explicit function $f : \bits^n \to \bits$, a computer mathematics package can easily calculate $\p_k(f; \epsilon)$ exactly, as a function of $k$ and $\epsilon$.  Furthermore, if $n$ is very small, a computer program can enumerate all antisymmetric left-monotone functions on $n$ bits.  We determined there are ``only'' 135 such functions for $n = 7$ and 2470 such functions for $n = 8$.  (The number jumps to 319124 for $n = 9$.)  Thus for particular small values of $n$ and $k$, we can completely solve the problem by comparing an explicit set of polynomials in $\eps$ on the range $(0, 1/2)$.

As an example, we give a nearly complete analysis of the case $n = 5$.  There are exactly 7 antisymmetric left-monotone functions on $n$ bits; they are $\maj{1}$, $\maj{3}$, $\maj{5}$, and four functions expressible as thresholds: $T_1 = Th(3, 1, 1, 1, 1; 4)$, $T_2 = Th(2, 1, 1, 1, 0 ; 3)$, $T_3 = Th(3, 2, 2, 1, 1; 5)$, and $T_4 = Th(4, 3, 2, 2, 2; 7)$, where $Th(a_1, \dots, a_5; \theta)$ is 1 iff $\sum_{i=1}^5 a_i x_i \geq \theta$.  Table~\ref{tab:1} shows $\p_k(f; \eps)$ for each of the functions.

For small values of $k$, we plotted these polynomials for $\eps \in (0,1/2)$.  This led to the following facts, which in principle could be proved by elementary analysis:

\begin{fact} \
\begin{itemize}
\item For $n = 5$, $2 \leq k \leq 9$, and any $\eps$, the best antisymmetric protocol is $\maj{1}$.
\item For $n = 5$, $k = 10, 11$, there exist $0 < \eps_k' < \eps_k'' < 1/2$ such that $\maj{3}$ is the best antisymmetric protocol for $\eps \in [\eps_k', \eps_k'']$, and $\maj{1}$ is the best antisymmetric protocol for all other $\eps$.
\item For $n = 5$, $k = 12$, there exist $0 < \eps_k' < \eps_k'' < \eps_k''' < 1/2$ such that $\maj{5}$ is the best antisymmetric protocol for $\eps \in [\eps_k', \eps_k'']$, $\maj{3}$ is the best antisymmetric protocol for $\eps \in [\eps_k'', \eps_k''']$, and $\maj{1}$ is the best antisymmetric protocol for all other $\eps$.
\end{itemize}
\end{fact}

The pattern for $k = 12$ appears to hold for all higher $k$, with $\maj{5}$ dominating more and more of the interval, as expected from Theorem~\ref{thm:majn}.

At this point we can prove two facts mentioned earlier:
\begin{proposition} \label{prop:majr}
There exist $k$, $\eps$, odd $n$, and odd $1 < r < n$ such that $\maj{r}$ is a superior protocol to both $\maj{1}$ and $\maj{n}$.
\end{proposition}
\begin{proof}
Substitute $k = 10$, $\eps = .26$ into Table~\ref{tab:1}.  By explicit calculation, $\p_{10}(\maj{1}; .26) \leq .0493$, $p_{10}(\maj{5}; .26) \leq .0488$, $p_{10}(\maj{3}; .26) \geq .0496.$
\end{proof}

\begin{proposition} \label{prop:unbalanced}
There exist $n$, $k$, and $f \in \B_k$ such that the probability all parties agree on 1 differs from the probability all parties agree on 0.
\end{proposition}
\begin{proof}
With $n = 5$, $k = 3$, and $f$ the left-monotone function with minterms $10010$ and $01101$, explicit calculation gives {\tiny $\frac12 - \frac{39}{16}\,\eps + 9\,\eps^2 - \frac{459}{16}\,\eps^3 + \frac{297}{4}\,\eps^4 - \frac{2331}{16}\,\eps^5 + \frac{3465}{16}\,\eps^6 - 234\,\eps^7 + 171\,\eps^8 - 75\,\eps^9 + 15\,\eps^{10}$} and {\tiny $\frac12 - \frac{39}{16}\,\eps + \frac{69}{8}\,\eps^2 - \frac{381}{16}\,\eps^3 + \frac{93}{2}\,\eps^4 - \frac{885}{16}\,\eps^5 + \frac{519}{16}\,\eps^6 + 6\,\eps^7 - 24\,\eps^8 + 15\,\eps^9 - 3\,\eps^{10}$} for the probabilities of agreement on 1 and 0, respectively.
\end{proof}

\begin{table}
\label{tab:1}
\begin{tabular}{|c|p{16cm}|}
\hline $f$ & $\p_k(f; \eps)$ \\\hline\hline

$\maj{1}$ & {\raggedright

\tiny  ${\epsilon}^{k}+ (1-\epsilon )^{k}$

\vspace*{-6ex}

} \\\hline

$T_1$ &  {\raggedright

\tiny  $1/16\, (-6\,{\epsilon}^{3}+5\,{\epsilon}^{4}-2\,{\epsilon}^{5}+4\,{
\epsilon}^{2} )^{k}+1/16\, (1+6\,{\epsilon}^{3}-5\,{\epsilon}^{4}
+2\,{\epsilon}^{5}-4\,{\epsilon}^{2} )^{k}+1/16\, (4\,\epsilon-10
\,{\epsilon}^{2}+10\,{\epsilon}^{3}-5\,{\epsilon}^{4}+2\,{\epsilon}^{5}
)^{k}+1/16\, (1-4\,\epsilon+10\,{\epsilon}^{2}-10\,{\epsilon}^{3}+5\,{
\epsilon}^{4}-2\,{\epsilon}^{5} )^{k}+1/4\, (\epsilon-{\epsilon}^{2
}+4\,{\epsilon}^{3}-5\,{\epsilon}^{4}+2\,{\epsilon}^{5} )^{k}+1/4\,
 (1-\epsilon+{\epsilon}^{2}-4\,{\epsilon}^{3}+5\,{\epsilon}^{4}-2\,{
\epsilon}^{5} )^{k}+1/4\, (1-2\,\epsilon+4\,{\epsilon}^{2}-6\,{
\epsilon}^{3}+5\,{\epsilon}^{4}-2\,{\epsilon}^{5} )^{k}+1/4\, (2
\,\epsilon-4\,{\epsilon}^{2}+6\,{\epsilon}^{3}-5\,{\epsilon}^{4}+2\,{\epsilon}^{
5} )^{k}+3/8\, (\epsilon+{\epsilon}^{2}-4\,{\epsilon}^{3}+5\,{
\epsilon}^{4}-2\,{\epsilon}^{5} )^{k}+3/8\, (1-\epsilon-{\epsilon}^
{2}+4\,{\epsilon}^{3}-5\,{\epsilon}^{4}+2\,{\epsilon}^{5} )^{k}$

\vspace*{-6ex}

} \\\hline

$T_2$ &  {\raggedright

\tiny  $1/8\, (-2\,{\epsilon}^{3}+3\,{\epsilon}^{2} )^{k}+1/8\, (1
+2\,{\epsilon}^{3}-3\,{\epsilon}^{2} )^{k}+1/8\, (3\,\epsilon-6\,
{\epsilon}^{2}+4\,{\epsilon}^{3} )^{k}+1/8\, (1-3\,\epsilon+6\,{
\epsilon}^{2}-4\,{\epsilon}^{3} )^{k}+3/8\,{\epsilon}^{k}+3/8\, (
1-\epsilon )^{k}+3/8\, (1-2\,\epsilon+3\,{\epsilon}^{2}-2\,{
\epsilon}^{3} )^{k}+3/8\, (2\,\epsilon-3\,{\epsilon}^{2}+2\,{
\epsilon}^{3} )^{k}$

\vspace*{-6ex}

} \\\hline

$\maj{3}$ &  {\raggedright

\tiny  $1/4\, (-2\,{\epsilon}^{3}+3\,{\epsilon}^{2} )^{k}+1/4\, (1
+2\,{\epsilon}^{3}-3\,{\epsilon}^{2} )^{k}+3/4\, (2\,\epsilon-3\,
{\epsilon}^{2}+2\,{\epsilon}^{3} )^{k}+3/4\, (1-2\,\epsilon+3\,{
\epsilon}^{2}-2\,{\epsilon}^{3} )^{k}$

\vspace*{-6ex}

} \\\hline

$T_3$ &  {\raggedright

\tiny $1/8\, (-6\,{\epsilon}^{3}+5\,{\epsilon}^{4}-2\,{\epsilon}^{5}+4\,{
\epsilon}^{2} )^{k}+1/8\, (1+6\,{\epsilon}^{3}-5\,{\epsilon}^{4}+
2\,{\epsilon}^{5}-4\,{\epsilon}^{2} )^{k}+1/16\, (\epsilon-{
\epsilon}^{2}+4\,{\epsilon}^{3}-5\,{\epsilon}^{4}+2\,{\epsilon}^{5} )^{k
}+1/16\, (1-\epsilon+{\epsilon}^{2}-4\,{\epsilon}^{3}+5\,{\epsilon}^{4}-2
\,{\epsilon}^{5} )^{k}+1/4\, (1-2\,\epsilon+4\,{\epsilon}^{2}-6\,
{\epsilon}^{3}+5\,{\epsilon}^{4}-2\,{\epsilon}^{5} )^{k}+1/4\, (2
\,\epsilon-4\,{\epsilon}^{2}+6\,{\epsilon}^{3}-5\,{\epsilon}^{4}+2\,{\epsilon}^{
5} )^{k}+1/8\, (1-\epsilon-{\epsilon}^{2}+4\,{\epsilon}^{3}-5\,{
\epsilon}^{4}+2\,{\epsilon}^{5} )^{k}+1/8\, (\epsilon+{\epsilon}^{2
}-4\,{\epsilon}^{3}+5\,{\epsilon}^{4}-2\,{\epsilon}^{5} )^{k}+3/16\,
 (8\,{\epsilon}^{3}-5\,{\epsilon}^{4}+2\,{\epsilon}^{5}+3\,\epsilon-7\,{
\epsilon}^{2} )^{k}+3/16\, (1-8\,{\epsilon}^{3}+5\,{\epsilon}^{4}
-2\,{\epsilon}^{5}-3\,\epsilon+7\,{\epsilon}^{2} )^{k}+3/16\, (2
\,{\epsilon}^{3}-5\,{\epsilon}^{4}+2\,{\epsilon}^{5}+2\,{\epsilon}^{2}+1-2\,
\epsilon )^{k}+3/16\, (-2\,{\epsilon}^{3}+5\,{\epsilon}^{4}-2\,{
\epsilon}^{5}-2\,{\epsilon}^{2}+2\,\epsilon )^{k}+1/16\, (2\,{
\epsilon}^{3}-5\,{\epsilon}^{4}+2\,{\epsilon}^{5}+2\,{\epsilon}^{2} )^{k
}+1/16\, (1-2\,{\epsilon}^{3}+5\,{\epsilon}^{4}-2\,{\epsilon}^{5}-2\,{
\epsilon}^{2} )^{k}$

\vspace*{-6ex}

} \\\hline

$T_4$ &  {\raggedright

\tiny $1/16\, ({\epsilon}^{2}+6\,{\epsilon}^{3}-10\,{\epsilon}^{4}+4\,{\epsilon}
^{5} )^{k}+1/16\, (1-{\epsilon}^{2}-6\,{\epsilon}^{3}+10\,{
\epsilon}^{4}-4\,{\epsilon}^{5} )^{k}+1/8\, (\epsilon+2\,{\epsilon}
^{2}-8\,{\epsilon}^{3}+10\,{\epsilon}^{4}-4\,{\epsilon}^{5} )^{k}+1/8
\, (1-\epsilon-2\,{\epsilon}^{2}+8\,{\epsilon}^{3}-10\,{\epsilon}^{4}+4\,
{\epsilon}^{5} )^{k}+1/16\, (1-2\,\epsilon+{\epsilon}^{2}+6\,{
\epsilon}^{3}-10\,{\epsilon}^{4}+4\,{\epsilon}^{5} )^{k}+1/16\, (
2\,\epsilon-{\epsilon}^{2}-6\,{\epsilon}^{3}+10\,{\epsilon}^{4}-4\,{\epsilon}^{5
} )^{k}+3/16\, (5\,{\epsilon}^{2}-10\,{\epsilon}^{3}+10\,{
\epsilon}^{4}-4\,{\epsilon}^{5} )^{k}+3/16\, (1-5\,{\epsilon}^{2}
+10\,{\epsilon}^{3}-10\,{\epsilon}^{4}+4\,{\epsilon}^{5} )^{k}+3/8\,
 (3\,\epsilon-8\,{\epsilon}^{2}+12\,{\epsilon}^{3}-10\,{\epsilon}^{4}+4\,
{\epsilon}^{5} )^{k}+3/8\, (1-3\,\epsilon+8\,{\epsilon}^{2}-12\,{
\epsilon}^{3}+10\,{\epsilon}^{4}-4\,{\epsilon}^{5} )^{k}+3/16\, (
1-2\,\epsilon+5\,{\epsilon}^{2}-10\,{\epsilon}^{3}+10\,{\epsilon}^{4}-4\,{
\epsilon}^{5} )^{k}+3/16\, (2\,\epsilon-5\,{\epsilon}^{2}+10\,{
\epsilon}^{3}-10\,{\epsilon}^{4}+4\,{\epsilon}^{5} )^{k}$

\vspace*{-6ex}

} \\\hline

$\maj{5}$ &   {\raggedright

 \tiny $1/16\, (10\,{\epsilon}^{3}-15\,{\epsilon}^{4}+6\,{\epsilon}^{5} )
^{k}+1/16\, (1-10\,{\epsilon}^{3}+15\,{\epsilon}^{4}-6\,{\epsilon}^{5}
 )^{k}+5/16\, (6\,{\epsilon}^{2}-14\,{\epsilon}^{3}+
15\,{\epsilon}^{4}-6\,{\epsilon}^{5} )^{k}+5/16\, (1
-6\,{\epsilon}^{2}+14\,{\epsilon}^{3}-15\,{\epsilon}^{4}+6\,{\epsilon}^{5}
 )^{k}+5/8\, (3\,\epsilon-9\,{\epsilon}^{2}+16\,{\epsilon}^{3}-15
\,{\epsilon}^{4}+6\,{\epsilon}^{5} )^{k}+5/8\, (1-3\,\epsilon+9\,
{\epsilon}^{2}-16\,{\epsilon}^{3}+15\,{\epsilon}^{4}-6\,{\epsilon}^{5} )
^{k}$

\vspace*{-6ex}

} \\\hline
\end{tabular}
\end{table}

We end with two open problems we've been led to consider:
\begin{open} Prove or disprove: For fixed $n$, $\eps$, and $2 \leq k \leq 9$, the best antisymmetric protocol is for all parties to use $\maj{1}$.
\end{open}
\begin{open} \label{conj:majlim} Prove or disprove: There is a universal constant $C < \infty$ such that for every $k, \eps$,
\[
\p_k(\maj{n^*}; \eps) \leq C \lim_{\substack{n\to \infty \\ n\text{ odd}}} \p(\maj{n},k,\eps),
\]
where $n^*$ is any odd number (presumably maximizing $\p_k(\maj{n^*} ; \eps)$).  I.e., the limiting value of $\p_k(\maj{n} ; \eps)$ is no worse than the success probability of the best majority, up to a constant factor.
\end{open}
The worst constant $C$ we know to be necessary in Open Problem~\ref{conj:majlim} is $\pi/2$, from the case $k = 2$, $\eps \to 1/2$.

If Conjecture O and Open Problem~\ref{conj:majlim} are both verified, then we can get tight (up to a constant) upper and lower bounds for the optimal value of $\p(f_1, \dots, f_k ; \eps)$ under antisymmetric protocols, for any $k$, $\eps$, and unrestricted $n$, by using Proposition~\ref{prop:integral}.\\

{\bf Acknowledgments:} We are most grateful to Oded Schramm for many helpful suggestions, ideas and proofs.
We are also grateful to Nati Srebro for helpful simulations and stimulating discussions.

\end{document}